\documentclass{amsart}
\usepackage{amssymb}
\usepackage{amsmath}
\usepackage{amsthm}
\newcommand{\pn}{\par\noindent}
\newcommand{\pmn}{\par\medskip\noindent}
\begin{document}
\title{On cubic polynomials with the cyclic Galois group}
\author{Yury Kochetkov}
\date{}
\email{yukochetkov@hse.ru,yuyukochetkov@gmail.com} \maketitle

\begin{abstract} A cubic Galois polynomial is a cubic polynomial with
rational coefficients that defines a cubic Galois field. Its
discriminant is a full square and its roots $x_1,x_2,x_3$
(enumerated in some order) are real. There exists (and only one)
quadratic polynomial $q$ with rational coefficients such that
$q(x_1)=x_2, q(x_2)=x_3, q(x_3)=x_1$. The polynomial $r=q(q)
\text{ mod } p$ cyclically permutes roots of $p$ in the opposite
order: $r(x_1)=x_3, r(x_3)=x_2, r(x_2)=x_1$. We prove that there
exist a unique Galois polynomial $p_1$ and a unique Galois
polynomial $p_2$ such that the polynomial $q$ cyclically permutes
roots of $p_1$ and the polynomial $r$ do the same with roots of
$p_2$. Polynomials $p$ and $p_1$ (and also $p$ and $p_2$) will be
called \emph{coupled}. Two polynomials are \emph{linear
equivalent}, if one of them is obtained from another by a linear
change of variable. By $C(p)$ we denote the class of polynomials,
linear equivalent to $p$. The coupling realizes a bijection
between classes $C(p)$ and $C(p_1)$ (and between classes $C(p)$
and $C(p_2)$). Classes $C(p)$ and $C(p_1)$ (and classes $C(p)$ and
$C(p_2)$) will be called \emph{adjacent}. We consider a graph: its
vertices --- are classes of the linear equivalency and two
vertices are connected by an edge, if the corresponded classes are
adjacent. Connected components of this graph will be called
\emph{superclasses}. In this work we give a description of
superclasses.\end{abstract}

\section{Coupled polynomials and classes of the linear equivalency}
\pn Let $p\in\mathbb{Q}[x]$ be an irreducible cubic polynomial (a
Galois polynomial) that defines a cubic Galois field. Roots
$x_1,x_2,x_3$ of such polynomial are real and its discriminant $D$
is a full square: $D=d^2$, $d\in\mathbb{Q}$. The Galois group of
$p$ is the cyclic group $A_3$ \cite{St}. \pmn {\bf Proposition 1.}
\emph{Let $x_1,x_2,x_3$ be roots of a Galois polynomial
$p=x^3+ax^2+bx+x$, enumerated in some order. There exists a unique
polynomial $q=\alpha x^2+\beta x+\gamma\in\mathbb{Q}[x]$ that
cyclically permutes roots of $p$: $q(x_1)=x_2$, $q(x_2)=x_3$,
$q(x_3)=x_1$.} \pmn \emph{Remark.} Let $K$ be the cubic Galois
field, generated by roots of $p$. The map $x\mapsto q(x)$ of $K$
into itself \emph{is not} an automorphism of $K$. \pmn
\emph{Remark.} The polynomial $q(q)\text{ mod }p$ permutes roots
of $p$ in the reverse order. \pmn \emph{Proof.} Let us consider
the linear system $$\left\{\begin{array}{l} \alpha x_1^2+\beta x_1+\gamma=x_2\\
\alpha x_2^2+\beta x_2+\gamma=x_3\\
\alpha x_3^2+\beta x_3+\gamma=x_1\end{array}\right.$$ The Cramer
formula gives us the solution of this system:
$$\alpha=\dfrac{\begin{vmatrix}
x_2&x_1&1\\x_3&x_2&1\\x_1&x_3&1\end{vmatrix}}
{\begin{vmatrix}x_1^2&x_1&1\\x_2^2&x_2&1\\x_3^2&x_3&1\end{vmatrix}}\,=
\dfrac{x_1^2+x_2^2+x_3^2-x_1x_2-x_1x_3-x_2x_3}{d}\,=
\frac{a^2-3b}{d}\,.$$ Analogously, we can find that
$$\beta=\frac{a^3+9c-7ab-d}{2d}\,,\,\gamma=\frac{a^2b+3ac-4b^2-ad}{2d}\,.$$
Thus,
$$q=\frac{a^2-3b}{d}\cdot x^2+\frac{a^3+9c-7ab-d}{2d}\cdot
x+\frac{a^2b+3ac-4b^2-ad}{2d}\,.\eqno(1)$$
Here the sign of $d$ is the sign of the number
$(x_1-x_2)(x_1-x_3)(x_2-x_3)$. If we choose another sign for $d$,
then we will have another solution of the above system. \qed \pmn
{\bf Example 1.} Let $p=x^3-3x+1$. The discriminant of $p$ is
$81$. The choice $d=9$ gives us the solution $q_1=-x^2-x+2$ and
the choice $d=-9$ --- the solution $q_2=x^2-2$. Obviously,
$q_1(q_1) \text{ mod p }=q_2$ and $q_2(q_2) \text{ mod } p=q_1$. \pmn
{\bf Corollary.}
\emph{The degree of $q$ is exactly 2.} \pmn \emph{Proof.} Let
$\alpha=0$, i.e. $3b=a^2$. Then $p=x^3+ax^3+\frac{a^2}{3}\cdot
x+c$ and $p'=3x^2+2ax+\frac{a^2}{3}=3\cdot(x+\frac a3)^2$. Thus,
the function $p(x)$ is nondecreasing, i.e. it has only one real
root and $p$ cannot be a Galois polynomial. \qed \pmn {\bf
Proposition 2.} \emph{Let $p$ be a cubic Galois polynomial and
polynomial $q$ cyclically permutes roots of $p$. Then $q$
cyclically permutes roots of another cubic Galois polynomial.}
\pmn \emph{Proof.} Let us consider the polynomial $s=(q(q(q)))-x$
of degree $8$. Each root of $p$ is a root of $s$, hence, $p$ is a
divisor of $s$. Each root of the polynomial $q-x$ is a root of
$s$, hence, $q-x$ is a divisor of $s$. Thus, there is the third
divisor of $s$ --- a polynomial $p_1$ of degree 3 with rational
coefficients. Let $x_1$ be a real root of $p_1$. Then it is a root
of $s$, i.e. $q(q(q(x_1)))=x_1$. Let $q(x_1)=x_2$ and
$q(x_2)=x_3$. As $q(q(q(x_2)))=x_2$ and $q(q(q(x_3)))=x_3$, then
$x_2$ and $x_3$ are roots of $p_1$. Thus, $q$ cyclically permutes
roots of $p_1$. But from (1) it follows that the discriminant of
$p_1$ is a full square. Thus, $p_1$ is a Galois polynomial. \qed
\pmn {\bf A continuation of Example 1.}  Here $p=x^3-3x+1$,
$q_1=-x^2-x+2$, $q_2=x^2-2$.
$$(q_1(q_1(q_1(x))))-x=(x^3-3x+1)(-x^2-2x+2)(x^3+2x^2-3x-5)$$ and
$$(q_2(q_2(q_2(x))))-x=(x^3-3x+1)(x^2-x-2)(x^3+x^2-2x-1).$$ The
discriminant of the polynomial $p_1=x^3+2x^2-3x-5$ is $169=13^2$
and the discriminant of the polynomial $p_2=x^3+x^2-2x-1$ is
$49=7^2$. Polynomials $p_1$ and $p_2$ define \emph{different}
Galois fields. \pmn \emph{Definition 1.} Two cubic Galois
polynomials $p$ and $r$ are called \emph{coupled}, if there exists
a quadratic polynomial $q$ that cyclically permutes roots of $p$
and roots of $r$. \pmn \emph{Remark.} As there are two polynomials
that cyclically permutes roots of a Galois polynomial $p$, then
$p$ is coupled with two Galois polynomials $p_1$ and $p_2$. \pmn
\emph{Definition 2.} Two polynomials are called linear equivalent,
if one of them is obtained from another by a linear change of
variable. The linear equivalency is an equivalency relation. The
set of polynomials, linear equivalent to a given polynomial $p$,
will be called the class of linear equivalency, generated by $p$,
and will be denoted $C(p)$. \pmn {\bf Proposition 3.} \emph{Let
$p$ and $r$ be coupled cubic Galois polynomial, $g(x)=p(\alpha
x+\beta)$ and $h(x)=r(\alpha x+\beta)$. Then $g$ and $h$ are
coupled cubic Galois polynomials.} \pmn \emph{Proof.} If the
polynomial $q$ cyclically permutes roots of $p$ and $r$, then the
polynomial $\frac{q(\alpha x+\beta)-\beta}{\alpha}$ cyclically
permutes roots of $g$ and $h$. \qed \pmn {\bf Corollary.}
\emph{The coupling is a bijection between $C(p)$ and $C(r)$.}

\section{Representatives of classes and characteristic numbers}
\pn \emph{Definition 3.}  Each class $C$ of linear equivalency
contains the unique polynomial of the form $x^3-ax-a$. This
polynomial will be called \emph{the representative} of the class
$C$. As the discriminant $D=a^2(4a-27)$ of this polynomial is a
full square, then $4a-27=k^2$. A rational number $k>0$ will be
called \emph{the characteristic number} of the class $C$. \pmn
{\bf Example 2.} Polynomial $x^3-27x-27$ is the representative of
the class $C(x^3-3x+1)$ and $9$ is the characteristic number of
this class. \pmn \emph{Remark.} Each cubic Galois field contains a
countable number of equivalency classes. For example, the field
generated by polynomial $x^3-3x+1$, contains equivalency classes
with representatives $x^3-tx-t$, where $t$ is any rational number
of the form
$$t=27\cdot\dfrac{(y^2+2187y+1594323)^3}{(y^3-4782969y-3486784401)^2}\,,
y\in\mathbb{Q}.$$ \pmn
{\bf Proposition 4.} \emph{Let $p=x^3-ax-a$, $a>0$, --- a cubic Galois
polynomial with discriminant
 $D=a^2k^2$ and let
$d=\sqrt D=ak$. Polynomials
$$q_1=\frac 3k\cdot x^2-\frac{k+9}{2k}\cdot x-\frac{2a}{k} \text{ and }
q_2=-\frac 3k\cdot x^2+\frac{9-k}{2k}+ \frac{2a}{k}\eqno(2)$$
induce cyclic permutations of roots of the polynomial $p$. Let
$p_1$ and $p_2$ be coupled polynomials. Polynomials
$$r_1=x^3-bx-b,\,\,b=\frac{27}{4}\cdot\dfrac{31k^2+108k+729}{(2k+27)^2}\,,
\text{ and }
r_2=x^3-cx-c,\,\,
c=\frac{27}{4}\cdot\dfrac{31k^2-108k+729}{(2k-27)^2}\eqno(3)$$ are
representatives of
classes $C(p_1)$ and $C(p_2)$.
The corresponding characteristic numbers are
$$k_1=\frac{27k}{2k+27} \text{ and }k_2=\frac{27k}{|2k-27|}\,.\eqno(4)$$}
\pmn \emph{Proof.}
Computation. \qed \pmn Thus, we have two maps in the set of
positive rational numbers $\mathbb{Q}_+$:
$$\varphi: k\mapsto \frac{27k}{2k+27} \text{ and } \psi: k\mapsto
\frac{27k}{|2k-27|}\,.\eqno(5)$$ \pmn
{\bf Proposition 5.} \emph{Maps $\varphi$ and $\psi$ have the
following properties:
\begin{enumerate} \item $\varphi(k)<k$, $\varphi(k)\in (0,\frac{27}{2})$;
\item iterations of
$\varphi(k)$ converge to zero; \item $\psi(\varphi(k))=k$;
$\varphi(\psi(k))=k$, if $k<\frac{27}{2}$; \item $\psi(k)>k$, if
$k<27$; $\psi(k)\in (\frac{27}{2},27)$, if $k>27$;
$\psi(\psi(k))=k$, if $k>\frac{27}{2}$\,; \item $\psi(27)=27$.
\end{enumerate}} \pmn\emph{Proof.} Only (2) needs a proof. We have,
$$\varphi(k)=\frac{27k}{2k+27}\,,\,\,\varphi(\varphi(k))=\frac{27k}{4k+27}\,,
\,\,
\varphi(\varphi(\varphi(k)))= \frac{27k}{6k+27}\,,\ldots$$ \qed
\pmn\emph{Remark.} Let $p$ be a cubic Galois polynomial and $p_1$
and $p_2$ be its coupled polynomials. Then $C(p_1)$ and $C(p_2)$
are different classes because their characteristic numbers are
different.

\section{Superclasses}
\pn \emph{Definition 4.} Two classes $C_1$ and $C_2$ will be
called \emph{adjacent} if there are coupled polynomials $p\in C_1$
and $r\in C_2$. \pmn \emph{Remark.} From Proposition 3 it follows
that if $C_1$ and $C_2$ are adjacent classes, then for each
element $g\in C_1$ there is a unique element $h\in C_2$, coupled
to $g$. \pmn \emph{Definition 5.} Let $G$ be a graph whose
vertices are classes of linear equivalency and two vertices are
connected by an edge, if corresponding classes are adjacent.
Connected components of $G$ will be called \emph{superclasses}.
\pmn {\bf Proposition 6.} \emph{Except two cases, each
superclass is generated by a positive rational number $k>27$ and
contains classes with characteristic numbers
$\{k,\psi(k),\varphi(k),\varphi(\psi(k)),\varphi(\varphi(k)),
\varphi(\varphi(\psi(k))),
\varphi(\varphi(\varphi(k))),\ldots\}$. Two exceptions are: a) the
superclass generated by $k=27$ (it contains classes with
characteristic numbers
$\{27,9,\frac{27}{5},\frac{27}{7},\ldots\}$); b) the superclass
generated by $k=\frac{27}{2}$ (it contains classes with
characteristic numbers
$\{\frac{27}{2},\frac{27}{4},\frac{27}{6},\ldots\}$).} \pmn
\emph{Remark.} Proposition 6 needs some clarification: a
superclass in our description is a set of characteristic numbers.
But it is possible, that some characteristic number in such set
corresponds to a class of reducible polynomials. For example,
number $k=270$ generates the superclass
$\{270,\frac{90}{7},\frac{270}{41},\frac{270}{61},\frac{10}{3},\ldots\}$.
Here the number $\frac{10}{3}$ corresponds to the class with
representative
$$x^3-\frac{343}{36}\cdot x-\frac{343}{36}=\left(x+\frac 73\right)
\left(x+\frac 76\right)\left(x-\frac 72\right).$$ It
must be noted that the coupled polynomial $1458x^3-7301x^2-6930x+49763$
is irreducible.

\vspace{5mm}

\end{document}